\documentclass{amsart}

\usepackage{amsmath}
\usepackage{amsfonts}
\usepackage{amssymb}
\usepackage{graphicx}
\usepackage{hyperref}
\usepackage{mathrsfs}
\usepackage{epstopdf}
\usepackage{amsmath,amsfonts,amsthm,enumerate,amscd,latexsym}
\usepackage{curves}
\usepackage{bbm}
\usepackage[mathscr]{eucal}
\usepackage{epsfig,epsf,xypic,epic}
\usepackage{caption}
\usepackage{subcaption}
\usepackage{tikz}
\usepackage{tikz-cd}
\usepackage[all]{xy}
\newtheorem{thm}{Theorem}[section]
\newtheorem{lem}[thm]{Lemma}

\newtheorem{prop}[thm]{Proposition}
\newtheorem{defn}[thm]{Definition}

\newtheorem{conj}[thm]{Conjecture}

\numberwithin{equation}{section}

\begin{document}

\title{From deformation theory to tropical geometry}
\author[K. Chan]{Kwokwai Chan}
\address{Department of Mathematics\\ The Chinese University of Hong Kong\\ Shatin\\ Hong Kong}
\email{kwchan@math.cuhk.edu.hk}

\date{\today}

\begin{abstract}
This is a write-up of the author's invited talk at the {\em Eighth International Congress of Chinese Mathematicians (ICCM)} held at Beijing in June 2019. We give a survey on the papers \cite{Chan-Leung-Ma17, Chan-Ma18} where the author and his collaborators Naichung Conan Leung and Ziming Nikolas Ma study how tropical objects arise from asymptotic analysis of the Maurer-Cartan equation for deformation of complex structures on a semi-flat Calabi-Yau manifold.
\end{abstract}

\maketitle

\tableofcontents

\section{Background}\label{sec:background}

This note is about a deep relationship between two apparently very different subjects: tropical geometry and deformation theory.


Recall that the complex structure $J$ on a compact complex manifold $M$ is an endomorphism of the tangent bundle $T_M$ of $M$ which squares to $-\text{Id}$. This induces and is equivalent to an eigenspace decomposition of the complexified tangent bundle $T_M\otimes_\mathbb{R} \mathbb{C} = T^{1,0}_M \oplus T^{0,1}_M$.
According to Kodaira-Spencer's classical theory \cite{Kodaira-Spencer-I-II, Kodaira-Spencer-III}, we deform $J$ by almost complex structures defined by elements $\Phi \in \text{Hom}(T^{0,1}_M, T^{1,0}_M) = \Omega^{0,1} (T^{1,0}_M)$. Such an almost complex structure is integrable if and only if $(\bar{\partial} + \Phi)^2 = 0$, which in turn is equivalent to the {\em Maurer-Cartan equation}
\begin{equation}\label{eqn:MC-eqn}
	\bar{\partial}\Phi + \frac{1}{2}[\Phi, \Phi] = 0
\end{equation}
associated to the {\em Kodaira-Spencer differential graded Lie algebra (DGLA)}
$$(\Omega^{0,\bullet}(M, T^{1,0}_M), \bar{\partial}, [\cdot,\cdot]).$$
It is a general philosophy that deformation problems are governed by DGLAs. In the above case, the space of infinitestimal deformations is given by the first cohomology group $H^1(M, T^{1,0}_M)$, while obstructions lie in the second cohomology group $H^2(M, T^{1,0}_M)$.

When $M$ is Calabi-Yau (i.e. $K_M \equiv \mathcal{O}_M$), one can enhance the DGLA to the {\em extended} Kodaira-Spencer complex
$$(\Omega^{0,\bullet}(M, \wedge^\bullet T^{1,0}_M), \bar{\partial}, \wedge, \Delta);$$
here $\Delta$ is the {\em Batalin-Vilkovisky (BV) operator} which corresponds to the operator $\partial$ on $\Omega^{0,\bullet}(M, \wedge^\bullet (T^{1,0}_M)^*) = \Omega^{\bullet,\bullet}(M)$ under the identification defined by contraction with the holomorphic volume form $\Omega$ on $M$. The BV operator $\Delta$ is not a derivation, but the discrepancy from being so is exactly measured by the Lie bracket:
$$[\alpha, \beta] = \Delta(\alpha \wedge \beta) - (\Delta \alpha)\wedge \beta - (-1)^{|\alpha|}\alpha \wedge (\Delta \beta).\footnote{This is also known as the Bogomolov-Tian-Todorov Lemma.}$$
This yields a so-called {\em differential graded Batalin–Vilkovisky (DGBV) algebra}, whose Maurer-Cartan equation \eqref{eqn:MC-eqn} is always solvable -- this is the famous unobstructedness result of Bogomolov-Tian-Todorov \cite{Bogomolov78, Tian86, Todorov89}.


On the other hand, tropical geometry is the study of algebraic geometry over the {\em tropical semi-field} $(\mathbb{T} := \mathbb{R} \cup \{+\infty\}, \oplus, \otimes)$, where the operations are defined by:
$$x \oplus y = \min \{x, y\},\qquad x \otimes y = x + y.$$
Tropical subvarieties are piecewise-linear objects, typical examples of which arise from tropical limits of classical subvarieties as follows: for $t > 0$, consider the map
$$\log_t: (\mathbb{C}^*)^n \to \mathbb{R}^n,\ (z_1, \dots, z_n) \mapsto (\log_t |z_1|, \dots, \log_t |z_n|).$$
Given an algebraic subvariety $X \subset (\mathbb{C}^*)^n$, the image $\mathcal{A}_t := \log_t (X)$ is called an {\em amoeba} of $X$. The tropical limit (or tropicalization) of $X$ is then given by $\Gamma = \lim_{t \to \infty} \mathcal{A}_t$, which is a tropical subvariety of $\mathbb{T}^n$. See Figure \ref{fig:amoeba} below.
\begin{figure}[ht]
	\includegraphics[scale=0.6]{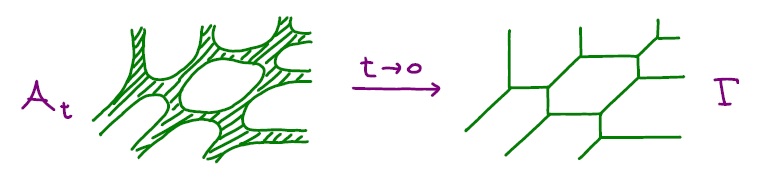}
	\caption{Amoeba of a holomorphic curve.}\label{fig:amoeba}
\end{figure}

Simplest examples of tropical varieties are tropical curves. For instance, a degree 2 polynomial in two variables is of the form
$$\min\{ a + 2x, b + x + y, c + 2y, d + x, e + y, f \},$$
whose zero set is defined as the set of points where the minimum is achieved by at least two entries, so it is given by an object as in the middle of Figure \ref{fig:tropical-curves}:
\begin{figure}[ht]
	\includegraphics[scale=0.55]{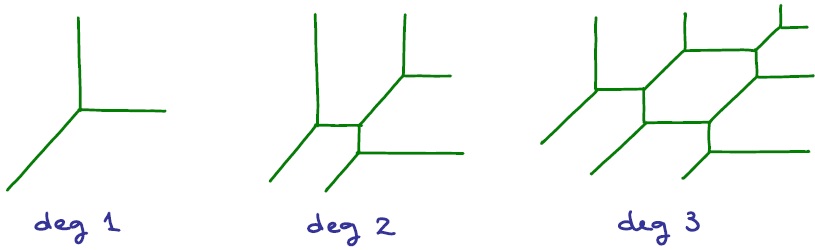}
	\caption{Tropical curves of low degrees.}\label{fig:tropical-curves}
\end{figure}

The first major application of tropical geometry was due to Mikhalkin \cite{Mikhalkin05}, who proved that counting holomorphic curves in toric surfaces is the same as counting tropical curves. This was later generalized to higher dimensions (using different methods) by Nishinou-Siebert \cite{Nishinou-Siebert06}. These results are proved by establishing a precise correspondence between holomorphic curves and tropical curves, from which one can clearly see a close relationship between Gromov-Witten theory and enumerative tropical geometry. However, it may not be at all obvious that there is a link between deformation theory and tropical geometry. 

\section{The link -- SYZ mirror symmetry}\label{sec:SYZ-Fukaya}

The key lies in mirror symmetry. More precisely, it was suggested in a 20 year old paper of Fukaya \cite{Fukaya05}, where he launched a spectacular program, expanding on the celebrated Strominger-Yau-Zaslow (SYZ) conjecture \cite{SYZ96}, which would give the ultimate explanation of the mechanism and miraculous power behind mirror symmetry. This section is a brief review of the story.


In 1996, Strominger, Yau and Zaslow \cite{SYZ96} made a ground-breaking proposal to explain mirror symmetry geometrically as a {\em T-duality}. In simple terms, what they asserted was that a mirror pair of Calabi-Yau manifolds should admit fiberwise dual (special) Lagrangian torus fibrations to the same base. In particular, this suggests the following mirror construction: given a (symplectic) Calabi-Yau manifold $X$ equipped with a Lagrangian torus fibration $\pi: X \to B$ admitting a Lagrangian section, the base manifold $B$ acquires an integral affine structure, and Duistermaat's global action-angle coordinates \cite{Duistermaat80} gives an identification
$$ X \cong T^*_B/\Lambda; $$
here $\Lambda$ is the natural lattice subbundle in $T^*_B$ generated by the coordinate $1$-forms $dx_1, \dots, dx_n$, where $\{x_1, \dots, x_n\}$ is a set of affine coordinates on $B$.
We can then define the {\em SYZ mirror} of $X$ as
$$ \check{X} := T_B / \Lambda^{\vee}, $$
where $\Lambda^{\vee} \subset T_B$ is the lattice dual to $\Lambda$ which is generated by the coordinate vector fields $\partial/\partial x_1, \dots, \partial/\partial x_n$. Since $B$ is an integral affine manifold, the quotient $\check{X}$ is naturally a complex manifold whose coordinates are given by exponentiation of complexification of the affine coordinates $\{x_1, \dots, x_n\}$. This demonstrates the mirror symmetry between $X$ and $\check{X}$ as a T-duality:
\begin{equation*}
\xymatrix{
	X = T^*_B/ \Lambda \ar[dr]_{\pi} & & T_B / \Lambda^\vee = \check{X} \ar[dl]^{\check{\pi}}\\
	& B &}
\end{equation*}

In general, however, the fibration $\pi: X \to B$ would have singular fibrations (or equivalently, the affine structure on the base $B$ would admit singularities). Thus the above construction can only be applied to the smooth locus $B_0 := B \setminus \Gamma$, where $\Gamma$ denotes the discriminant locus of the fibration $\pi$, which gives at best an approximation of the true picture:
\begin{equation}\label{eqn:dual-fibration-smooth-locus}
\xymatrix{
	X \supset X_0 = T^*_{B_0}/ \Lambda \ar[dr]_{\pi} & & T_{B_0} / \Lambda^\vee = \check{X}_0 \ar[dl]^{\check{\pi}}\\
	& B_0 &}
\end{equation}

Of course $\check{X}_0$ is not the correct mirror because singular fibers of $\pi$ have all been removed from $X$ and information is lost. We expect the correct mirror $\check{X}$ to be given by a (partial) compactification of $\check{X}_0$. But then a problem arises: the natural complex structure $\check{J}_0$ on $\check{X}_0$ can {\em never} be extended to {\em any} partial compactification of $\check{X}_0$ due to nontrivial monodromy of the affine structure on $B$ around the discrminant locus $\Gamma$. This leads to the most crucial idea in the SYZ proposal \cite{SYZ96}: one should use {\em quantum corrections} coming from holomorphic disks in $X$ with boundaries on Lagrangian torus fibers of $\pi$ to correct or deform $\check{J}_0$ so that it becomes extendable. This is the so-called {\em reconstruction problem} in mirror symmetry.

\begin{figure}[ht]
	\includegraphics[scale=0.45]{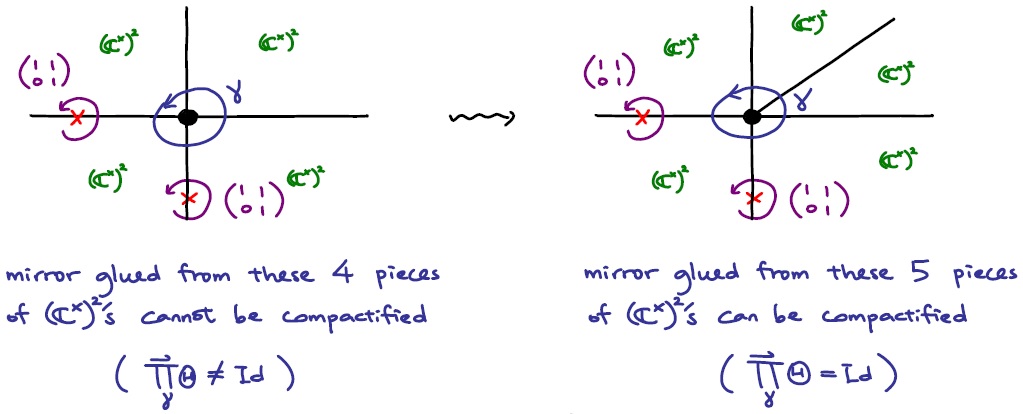}
	\caption{Corrected gluing of the mirror using a consistent scattering diagram.}\label{fig:gluing-scattering}
\end{figure}

This problem was solved in the 2-dimensional case (over non-Archimedean fields) by Kontsevich-Soibelman \cite{Kontsevich-Soibelman06} and in general dimensions (over $\mathbb{C}$) by Gross-Siebert \cite{Gross-Siebert-reconstruction}. In these fundamental works, a class of tropical objects called {\em scattering diagrams} was used to describe the quantum corrections used to modify the coordinate changes (or gluings) in $\check{X}_0$ along walls in the base (see Figure \ref{fig:gluing-scattering}). In other words, they tackled the reconstruction problem via a \v{C}ech approach. This is made possible by the following key lemma due to Kontsevich-Soibelman:
\begin{lem}[Kontsevich-Soibelman \cite{Kontsevich-Soibelman06}]
	Any scattering diagram $\mathscr{D}_0$ can be completed (by adding rays) to a {\em consistent} scattering diagram $\mathscr{D}$, i.e.
	$$\Theta_\gamma := \prod^{\rightarrow}_\gamma \Theta = \text{Id}$$
	for any loop $\gamma$ around any singular point of $\mathscr{D}$; here $\prod^{\rightarrow}$ is the path-ordered product of the automorphisms $\Theta$ associated to the rays of $\mathscr{D}$ which intersect $\gamma$.
\end{lem}

Both Kontsevich-Soibelman \cite{Kontsevich-Soibelman06} and Gross-Siebert \cite{Gross-Siebert-reconstruction} were actually motivated by Fukaya's proposal \cite{Fukaya05}.\footnote{Similar ideas, but using rigid analytic geometry, instead of asymptotic analysis, appeared in an even earlier work of Kontsevich-Soibelman \cite{Kontsevich-Soibelman01}.}
In his differential-geometric andf more transcendental approach \cite{Fukaya05}, Fukaya considered the Kodaira-Spencer DGLA
$$(\Omega^{0,\bullet}(\check{X}_0, T^{1,0}_{\check{X}_0}), \bar{\partial}, [\cdot,\cdot])$$
on the semi-flat SYZ mirror $\check{X}_0$ and the associated Maurer-Cartan equation \eqref{eqn:MC-eqn}. 
He packaged the quantum corrections into a $(0,1)$-form $\Phi$ with values in the holomorphic tangent bundle $T^{1,0}_{\check{X}_0}$ which solves the Maurer-Cartan equation, or equivalently, as a deformation in the classical Kodaira-Spencer theory.
And then he studied what happens when Maurer-Cartan solutions are expanded into Fourier series along the Lagrangian torus fibers of $\check{\pi}: \check{X}_0 \to B_0$ (the dual of $\pi: X_0 \to B_0$ as shown in the diagram \eqref{eqn:dual-fibration-smooth-locus}). 
In this way, Fukaya related holomorphic disks in $X$ with boundaries on Lagrangian torus fibers of $\pi: X \to B$ with deformations of the complex structure $\check{J}_0$ on the semi-flat SYZ mirror $\check{X}_0$, via Morse theory on the base manifold $B$, and made a series of spectacular conjectures including the following:
\begin{conj}[Fukaya \cite{Fukaya05}]
	Near a large complex structure limit, the Fourier modes of the Maurer-Cartan solutions are supported near gradient flow trees of the area functional (which is a multi-valued Morse function) on $B$, and these gradient flow trees in $B$ are adiabatic limits of holomorphic disks in $X$ with boundaries on Lagrangian torus fibers of $\pi: X \to B$.
\end{conj}
\begin{figure}[ht]
	\includegraphics[scale=0.55]{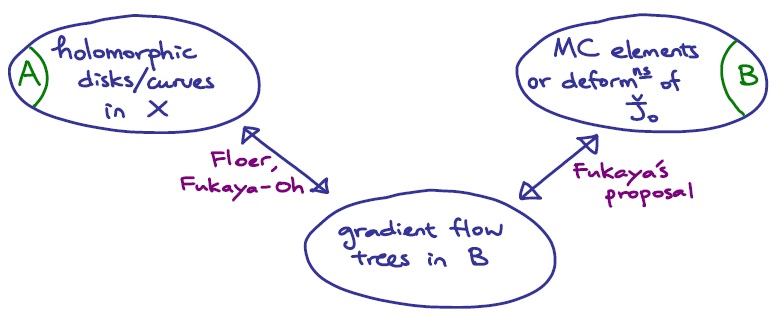}
	\caption{}\label{fig:fukaya-proposal-I}
\end{figure}
One can imagine that holomorphic disks can glue to produce holomorphic curves in $X$, so Fukaya's conjectures beautifully explain why Gromov-Witten theory of $X$ is encoded in the deformation theory of the mirror $\check{X}$ and why mirror symmetry can be applied to make enumerative predictions.
This is a much more transparent description of the picture depicted by the SYZ conjecture \cite{SYZ96}. Unfortunately, the arguments in \cite{Fukaya05} were only heuristical and the analysis involved to make them precise seemed intractable at that time.

To reduce the analytical difficulties, we observe that tropical objects, perhaps because of their linear nature, are much easier to deal with than Morse-theoretical objects.\footnote{In retrospect, this was a lesson we learned from the Gross-Siebert program \cite{Gross-Siebert03, Gross-Siebert-logI, Gross-Siebert-logII, Gross-Siebert-reconstruction}.} Moreover, we lose no information by traveling to the tropical world because tropical geometry on $B$ and Morse theory on $B$ are equivalent via a Legendre transform. For instance, gradient flow trees in $B$ correspond to tropical trees in $B$.
In view of this, it seems conceivable to modify Fukaya's proposal as a relationship between deformation theory on $\check{X}_0$ and tropical geometry on $B$, as illustrated by Figure \ref{fig:fukaya-proposal-II} below.
\begin{figure}[ht]
	\includegraphics[scale=0.88]{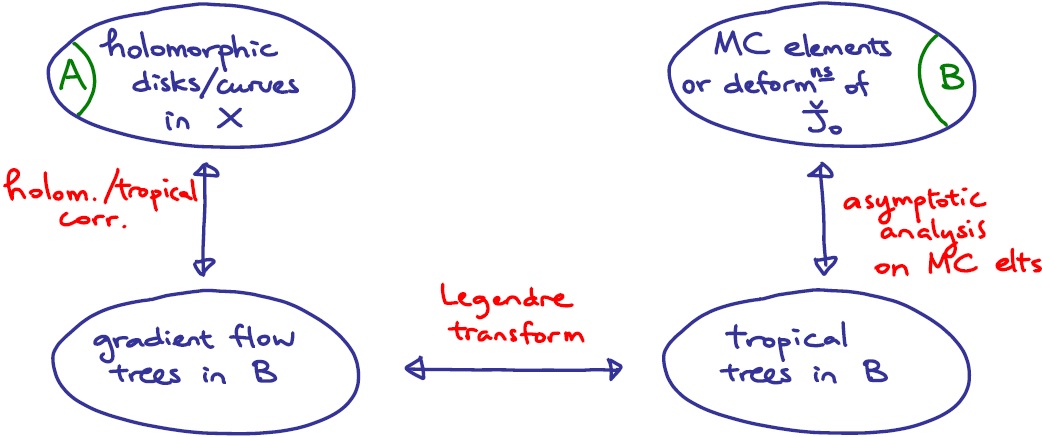}
	\caption{}\label{fig:fukaya-proposal-II}
\end{figure}

\section{Scattering diagrams from Maurer-Cartan elements}\label{sec:scattering-from-MC}

In \cite{Chan-Leung-Ma17}, Leung, Ma and the author made the first attempt to realize Fukaya's (modified) proposal. We established a precise relation between consistent scattering diagrams in $B$ and Maurer-Cartan solutions of the Kodaira-Spencer DGLA of the semi-flat Calabi-Yau manifold $\check{X}_0$. Since consistent scattering diagrams can be thought geometrically as tropical limits of loci (in the SYZ base $B$) of Lagrangian torus fibers of $\pi: X \to B$ which bound (Maslov index $0$) holomorphic disks in $X$, this describes a local model for how SYZ mirror symmetry works.

To explain our results, we first recall the definition of a scattering diagram. 
To simplify the exposition, we restrict ourselves to the $2$-dimensional case in the rest of this note. While the 2d case is sufficient for illustrating the relationship between deformation theory and tropical geometry, we shall emphasize that many of the notions and results below actually work in higher dimensions as well. We will mainly be following the notations and definitions in \cite{Gross10, GPS10}; see also \cite{Mikhalkin05, Nishinou-Siebert06, Nishinou12}.

We fix, once and for all, a lattice $M \cong \mathbb{Z}^2$ with basis $e_1 = (1,0)$ and $e_2 = (0,1)$. Let $N = \text{Hom}(M, \mathbb{Z})$ be the dual lattice, and $M_\mathbb{R} = M \otimes_\mathbb{Z} \mathbb{R}$ and $N_\mathbb{R} = N \otimes_\mathbb{Z} \mathbb{R}$ be the associated real vector spaces. 
For $m = (a, b) \in M$, we denote by $z^m = z_1^a z_2^b \in \mathbb{C}[M] \cong \mathbb{C}[z_1^\pm, z_2^\pm]$ the corresponding monomial.

We take $B_0 = M_\mathbb{R} \cong \mathbb{R}$ and identify it with $N_\mathbb{R}$ using the flat metric on $\mathbb{R}^2$. We also equip $B_0$ with the positive orientation. Accordingly, we have $\check{X}_0 \cong (\mathbb{C}^*)^2$ and the fibration $\check{\pi}: \check{X}_0 \to B_0$ is nothing but the log map $\log: (\mathbb{C}^*)^2 \to \mathbb{R}^2,\ (z_1, z_2) \mapsto (\log |z_1|, \log |z_2|)$.

Consider $\mathfrak{g} := (\mathbb{C}[M] \widehat{\otimes}_\mathbb{C} \mathbf{m}) \otimes_\mathbb{Z} N$, where $\mathbf{m}$ is the maximal ideal $(t)$ in the power series ring $\mathbb{C}[[t]]$, equipped with the Lie bracket
$$[z^m \partial_n, z^{m'}\partial_{n'}] := z^{m + m'} \partial_{\langle m',n \rangle n' - \langle m,n' \rangle n},$$
and the Lie subalgebra 
$$\mathfrak{h} := \bigoplus_{m \in M \setminus \{0\}} ((\mathbb{C}\cdot z^m) \widehat{\otimes}_\mathbb{C}\mathbf{m}) \otimes_\mathbb{Z} m^\perp.$$
Via exponentiation, this defines the so-called {\em tropical vertex group} \cite{GPS10}
$$\mathbb{V} := \exp(\mathfrak{h}).$$
Elements of $\mathbb{V}$ should be viewed as automorphisms of (formal families) of $(\mathbb{C}^*)^2$, which are used to correct the gluings in the reconstruction problem; see Figure \ref{fig:gluing-scattering}. We always write an element $\Theta \in \mathbb{V}$ as $\Theta = \exp(f\partial_n)$, where 
$$f = \sum_{j,k \geq 1} a_{jk} z^{km} t^{j} \in (\mathbb{C}[z^m]\cdot z^m)\widehat{\otimes}_\mathbb{C} \mathbf{m}$$
and $n \in m^\perp \subset N$.
\begin{defn}\label{def:wall}
A {\em wall} is a triple $(m, P, \Theta)$, where
\begin{itemize}
	\item
	$m = (a, b) \in M \setminus \{0\}$,
	
	\item
	$P \subset M_\mathbb{R}$ is either a line of the form $P = p + \mathbb{R} m$ or a ray of the form $P = p + \mathbb{R}_{\geq 0} m$, and
	
	\item
	$\Theta = \exp(f\partial_n) \in \mathbb{V}$ such that $n \in m^\perp \setminus \{0\}$ is the unique primitive vector so that $\{m ,n\}$ defines the positive orientation on $M_\mathbb{R} \cong \mathbb{R}^2$ (after identifying $N_\mathbb{R}$ with $M_\mathbb{R}$ using the flat metric on $\mathbb{R}^2$).
\end{itemize}

We call $\Theta$ the {\em wall-crossing factor} associated to the wall $\mathbf{w}$.

\end{defn}
\begin{defn}\label{def:scattering-diagram}
	$\text{ }$
	\begin{itemize}
		\item A {\em scattering diagram} $\mathscr{D}$ is a set of walls $\left\{ ( m_\alpha, P_\alpha, \Theta_\alpha) \right\}_{\alpha }$ such that there are only finitely many $\alpha$'s with $\Theta_\alpha \neq id $ $(\text{mod $\mathbf{m}^N$})$ for every $N \in \mathbb{Z}_{>0}$. 
		\item We define the {\em support} of a scattering diagram $\mathscr{D}$ to be $
		\text{Supp}(\mathscr{D}) := \bigcup_{\mathbf{w} \in \mathscr{D}} P_{\mathbf{w}}$, and the {\em singular set} of $\mathscr{D}$ to be $
		\text{Sing}(\mathscr{D}) := \bigcup_{\mathbf{w} \in \mathscr{D}} \partial P_{\mathbf{w}} \cup \bigcup_{\mathbf{w}_1\pitchfork \mathbf{w}_2} P_{\mathbf{w}_1} \cap P_{\mathbf{w}_2}$,
		where $\mathbf{w}_1 \pitchfork \mathbf{w}_2$ means transversally intersecting walls.
		\item A scattering diagram $\mathscr{D} = \left\{ ( m_\alpha, P_\alpha, \Theta_\alpha) \right\}_{\alpha }$ is said to be {\em consistent} if for any loop $\gamma$ around a singular point of $\mathscr{D}$, the path-ordered product along $\gamma$ (i.e. product of elements $\Theta_\alpha$ associated to walls which intersect the loop $\gamma$ in the order determined by the orientation on $\gamma$) is the identity:
		$$\Theta_\gamma := \prod^{\rightarrow}_\gamma \Theta_\alpha = \text{Id}.$$
	\end{itemize}
\end{defn}

\begin{figure}[ht]
	\includegraphics[scale=0.45]{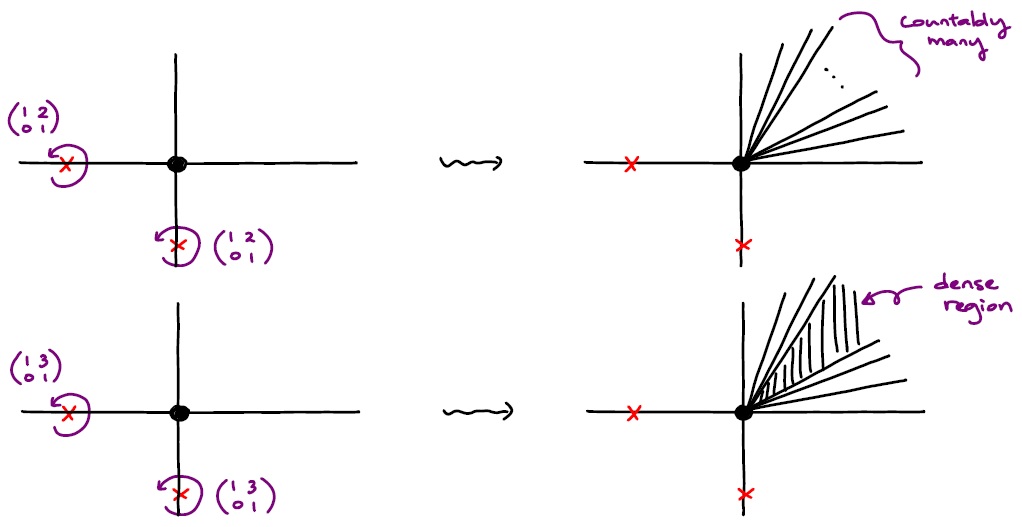}
	\caption{Examples of consistent scattering diagrams from completing inconsistent two-wall diagrams.}\label{fig:scattering-examples}
\end{figure}

The simplest nontrivial example of a consistent scattering diagram is the {\em pentagon diagram} shown on the RHS of Figure \ref{fig:gluing-scattering}. In general, the combinatorics of such diagrams can be very complicated; Figure \ref{fig:scattering-examples} shows a couple more examples.


To see how these are related to deformation theory, we first take the Fourier expansion of the Kodaira-Spencer DGLA $(\Omega^{0,\bullet}(\check{X}_0, T^{1,0}_{\check{X}_0}), \bar{\partial}, [\cdot,\cdot])$ along the Lagrangian torus fibers of $\check{\pi}: \check{X}_0 \to B_0$. This produces a DGLA
$$\left(L^\bullet = \bigoplus_{i \in \mathbb{Z}_{\geq 0}} L^i, \bar{\partial}, [\cdot,\cdot]\right)$$
over the integral affine manifold $B_0$.

We start with the simplest (and trivial) consistent scattering diagram, namely, the diagram consisting of just one wall $\mathbf{w} = (m, P, \Theta)$, where $P = p + \mathbb{R} m$ is a line. In \cite{Chan-Leung-Ma17}, we showed that an explicit solution (or ansatz) of the Maurer-Cartan equation \eqref{eqn:MC-eqn} associated to $\mathbf{w}$ can be written down. Moreover, asymptotic analysis of the Fourier modes of this Maurer-Cartan solution can recover the wall-crossing factor $\Theta$:
\begin{prop}[\S4 in \cite{Chan-Leung-Ma17}]
	Let $\mathbf{w} = (m, P, \Theta)$ be a wall supported on a line. Then we have the following statements.
	\begin{enumerate}
		\item There exists a Maurer-Cartan solution of the form
		$$\Xi_\mathbf{w} = -\delta_m\cdot \log(\Theta) \in L^1,$$
		where $\delta_m$ is the bump function (i.e. smoothing of a delta function) along the direction orthogonal to $m$ (see Figure \ref{fig:one-wall-ansatz} below).
		\item The Maurer-Cartan solution $\Xi_\mathbf{w}$ is gauge equivalent to $0$, and after choosing a suitable gauge fixing condition, there is a unique element $\varphi_\mathbf{w} \in L^0$ such that $e^{\varphi_\mathbf{w}} \ast 0 = \Xi_\mathbf{w}$.
		\item We have the asymptotic expansion
		$$\varphi_\mathbf{w} = (\varphi_\mathbf{w})_0 + O(\hbar)$$
		as $\hbar \to 0$; furthermore, the Fourier series of the leading term $(\varphi_\mathbf{w})_0$ is of the form
		\begin{equation*}
			\mathcal{F}((\varphi_\mathbf{w})_0) = \left\{
				\begin{array}{ll}
					\log(\Theta) & \text{on $H_+$,}\\
					0 & \text{on $H_-$.}
				\end{array}\right.
		\end{equation*}
	\end{enumerate}
\end{prop}
\begin{figure}[ht]
	\includegraphics[scale=0.55]{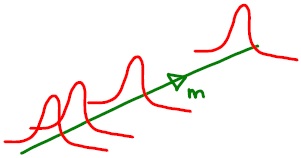}
	\caption{Ansatz for solving the MC equation in the one wall case.}\label{fig:one-wall-ansatz}
\end{figure}
Here, $\hbar$ measures the size of fibers of the fibration $\check{\pi}: \check{X}_0 \to B_0$, and $\hbar \to 0$ is a {\em large complex structure limit}; $H_\pm$ are half-spaces in $B_0 = M_\mathbb{R}$ defined by $P$.


Next we study the much more nontrivial and important case of two transversally intersecting walls: $\mathbf{w}_1 = (m_1, P_1, \Theta_1)$, $\mathbf{w}_2 = (m_2, P_2, \Theta_2)$, where $P_1 = p + \mathbb{R} m_1$ and $P_2 = p + \mathbb{R} m_2$ are lines intersecting at the point $p \in B_0$. Now the superposition
$$\Pi := \Xi_{\mathbf{w}_1} + \Xi_{\mathbf{w}_2} \in L^1$$
does not solve the Maurer-Cartan equation \eqref{eqn:MC-eqn} around the intersection point $p$, although each summand itself is a solution.
\begin{figure}[ht]
	\includegraphics[scale=0.55]{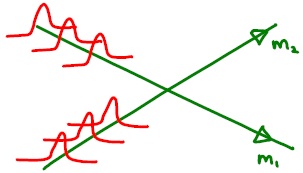}
	\caption{Two walls intersecting.}\label{fig:two-walls-intersecting}
\end{figure}

To solve the Maurer-Cartan equation \eqref{eqn:MC-eqn} by correcting $\Pi$, we apply {\em Kuranishi's method} \cite{Kuranishi65}. The key ingredient is a {\em propagator} (or {\em gauge fixing}) $H: L^\bullet \to L^\bullet$, i.e. a homotopy retract 
\begin{equation*}
	\begin{tikzcd}
		H^*(L^\bullet) \arrow[r,shift left=.5ex,"\iota"]
		&
		L^\bullet \arrow[l,shift left=.5ex,"p"]
		\circlearrowleft H
	\end{tikzcd}
\end{equation*}
such that
\begin{align*}
\text{Id} - p \circ \iota  = 0,\quad
\text{Id} - \iota \circ p  = dH + Hd.
\end{align*}
A solution of \eqref{eqn:MC-eqn} is then given by the formula
\begin{equation}\label{eqn:Kuranishi-formula}
\Phi := \Pi - \frac{1}{2} H [\Phi, \Phi],
\end{equation}
upon checking $p[\Phi, \Phi] = 0$.
An interesting and important fact is that the RHS of Kuranishi's formula \eqref{eqn:Kuranishi-formula} can be expressed as a sum over trivalent trees:
\begin{equation}\label{eqn:sum-over-tree}
\Phi = \sum_{T:\text{trivalent trees}} W_T,
\end{equation}
where each summand $W_T$ is computed by aligning the input $\Pi$ at each of the leaves of the tree $T$, taking the Lie bracket $[\cdot,\cdot]$ of two incoming terms at a trivalent vertex, and also acting by the propagator $H$ on any internal edge; see Figure \ref{fig:trivalent-tree} below.
\begin{figure}[ht]
	\includegraphics[scale=0.55]{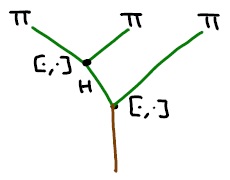}
	\caption{A trivalent tree which can appear in the summation formula for $\Phi$.}\label{fig:trivalent-tree}
\end{figure}

Indeed, it is well-known that the Maurer-Cartan equation associated to any DGLA can be solved by such sum-over-trees formulas. This is also a commonly used technique in solving the classical or quantum master equation in quantum field theory.

The upshot is that now the Maurer-Cartan solution \eqref{eqn:Kuranishi-formula} admits a Fourier expansion which naturally gives rise to a consistent scattering diagram that completes the two-wall diagram:
\begin{thm}[Theorem 1.5 in \cite{Chan-Leung-Ma17}]\label{thm:Chan-Leung-Ma}
	Let $\mathbf{w}_1 = (m_1, P_1, \Theta_1)$, $\mathbf{w}_2 = (m_2, P_2, \Theta_2)$ be two walls intersecting transversally.  Setting $\Pi := \Xi_{\mathbf{w}_1} + \Xi_{\mathbf{w}_2} \in L^1$ and $\Phi := \Pi - \frac{1}{2} H [\Phi, \Phi]$, then we have the following statements.
	\begin{enumerate}
		\item The Fourier expansion of $\Phi$ is of the form
		$$\Phi = \Pi + \sum_{m \in \left(\mathbb{Z}^2_{>0}\right)_{\text{prim}}} \Phi_{m},$$
		where each summand $\Phi_m$ is a Maurer-Cartan solution supported near the ray $P_m := p + \mathbb{R}_{\geq 0}m$.\footnote{To make precise the phrase ``supported near'', we introduced the notion of {\em asymptotic support} in \cite[Definition 4.19]{Chan-Leung-Ma17}.}
		\item After choosing a suitable gauge fixing condition, there exists, for each $m \in \left(\mathbb{Z}^2_{>0}\right)_{\text{prim}}$, a unique element  $\varphi_m \in L^0$ such that
				\begin{enumerate}
					\item $e^{\varphi_m} \ast 0 = \Phi_m$,
					\item $\varphi_m = (\varphi_m)_0 + O(\hbar)$ asymptotically as $\hbar \to 0$, and
					\item there exists an element $\Theta_m \in \mathbb{V}$ so that the Fourier series of $(\varphi_m)_0$ is precisely given by
					\begin{equation*}
						\mathcal{F}((\varphi_m)_0) = \left\{
						\begin{array}{ll}
							\log(\Theta_m) & \text{on $H_{m,+}$,}\\
							0 & \text{on $H_{m,-}$.}
						\end{array}\right.
					\end{equation*}
				\end{enumerate}
		From these, we can associate a scattering diagram $\mathscr{D}(\Phi)$ to the Maurer-Cartan solution $\Phi$.
		\item The scattering diagram $\mathscr{D}(\Phi)$ is consistent.
	\end{enumerate}
\end{thm}
\begin{figure}[ht]
	\includegraphics[scale=0.55]{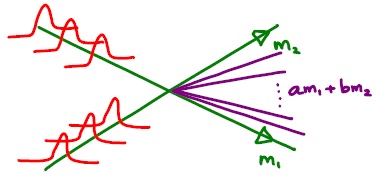}
	\caption{Completion of the two-wall diagram by solving the MC equation}\label{fig:two-walls-completed}
\end{figure}

A scattering diagram can be viewed as a number of tropical disks stacked together (cf. the notion of standard scattering diagram in Definition 1.10 and the subsequent discussion in \cite[\S 1]{GPS10}). From the proof of Theorem \ref{thm:Chan-Leung-Ma} in \cite{Chan-Leung-Ma17}, one can actually see that each trivalent tree in \eqref{eqn:sum-over-tree} gives rise to exactly one such tropical disk.

\section{Tropical counting from Maurer-Cartan elements}\label{sec:tropical-from-MC}

An application of the ideas in the previous section can reveal further connections between tropical geometry and deformation theory. In \cite{Chan-Ma18}, Ma and the author studied mirror symmetry for toric Fano surfaces via Fukaya's program. Here we briefly review what we have obtained.

Recall that the mirror of a toric Fano surface $X = X_\Sigma$ is given by a {\em Landau-Ginzburg model} $(\check{X}, W)$ where $\check{X} = \check{X}_0 = (\mathbb{C}^*)^2$ and $W: \check{X} \to \mathbb{C}$ is a Laurent polynomial called the {\em Hori-Vafa superpotential} \cite{Hori-Vafa00}. In an important work \cite{Cho-Oh06}, Cho-Oh showed that $W$ coincides with the {\em Lagrangian Floer superpotential} of $X$, defined via Lagrangian Floer theory by Fukaya-Oh-Ohta-Ono \cite{FOOO-book, FOOO-toricI, FOOO-toricII, FOOO-toricIII}. More previsely, they showed that the coefficients of $W$ are exactly counts of Maslov index $2$ holomorphic disks in $X$ with boundaries on Lagrangian torus fibers of the moment map $\pi: X \to \Delta$, where $\Delta \subset M_\mathbb{R}$ is the moment polytope.

\begin{figure}[ht]
	\includegraphics[scale=0.65]{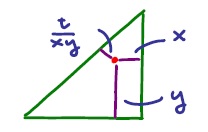}
	\caption{Moment polytope of $\mathbb{P}^2$ and the MI 2 tropical disks.}\label{fig:tropical-disks-P2}
\end{figure}

A prototypical example is given by mirror symmetry of the projective plane $X = \mathbb{P}^2$. The mirror superpotential is explicitly given by the Laurent polynomial $W = x + y + \frac{t}{xy}$ on $\check{X} = (\mathbb{C}^*)^2$. One can directly see that the monomial terms in $W$ are in a one-to-one correspondence with Maslov index $2$ disks in $X$ bounded by moment map Lagrangian torus fibers; see Figure \ref{fig:tropical-disks-P2} for the tropicalization of such disks.

In the toric case, the simplest nontrivial mirror statement is the ring isomorphism $QH^*(X) \cong Jac(W)$, where $QH^*(X)$ is the {\em small} quantum cohomology ring of $X$ and $Jac(W) := \mathbb{C}[x^{\pm}, y^{\pm}] / (x\partial W/\partial x, y\partial W/\partial y)$ is the Jacobian ring of $W$ \cite{Batyrev93, Givental98} (see also \cite{Chan-Leung10a}). To extend this to the {\em big} quantum cohomology, we need to perturb the superpotential $W$. There are two methods: one is by bulk deformations, due to Fukaya-Oh-Ohta-Ono \cite{FOOO-toricII}; another is by counts of tropical disks, due to Gross \cite{Gross10}, but only in the case of $\mathbb{P}^2$. We will recall the latter approach briefly below.

Following \cite{Gross10}, we consider a graph $\Gamma$ without bivalent vertices but with some unbounded edges. Denote by $\Gamma^{[0]}$ and $\Gamma^{[1]}$ the set of vertices and the set of edges of $\Gamma$ respectively; also denote by $\Gamma^{[1]}_\infty \subset \Gamma^{[1]}$ the set of unbounded edges. Let $X_\Sigma$ be a toric Fano surface. Then a {\em $d$-pointed marked tropical curve in $X_\Sigma$} is a map $h: \Gamma \to M_\mathbb{R}$ together with a weight function $w: \Gamma^{[0]} \to \mathbb{Z}_{\geq 0}$ and a marking $\{p_1, \dots, p_d\} \hookrightarrow \Gamma^{[1]}_\infty$ satisfying the following conditions:
\begin{itemize}
	\item $w(E) = 0$ if and only if $E = E_{p_i}$ for some $i$.
	\item For $i = 1, \dots, d$, $h|_{E_{p_i}}$ is constant, while for each $E \in \Gamma^{[1]} \setminus \{E_{p_1}, \dots, E_{p_d}\}$, $h|_E$ is a line segment of rational slope.
	\item At each vertex $V \in \Gamma^{[0]}$, the {\em balancing condition} is satisfied, meaning that, if $E_1, \dots, E_\ell \in \Gamma^{[1]}$ are the edges adjacent to $V$, and $m_i \in M$ denotes the primitive vector tangent to $h(E_i)$ and pointing away from $V$, then we have
	$$\sum_{i = 1}^\ell w(E_i)m_i = 0.$$
	\item For each unbounded edge $E \in \Gamma^{[1]}_\infty$, $h(E)$ is either a point or an affine translate of some $1$-dimensional cone $\rho$ in the fan $\Sigma$.
\end{itemize}

To define tropical disks, we consider a graph of the form $\Gamma' = \Gamma \cup \{V_{\text{out}}\}$, where $\Gamma$ is as above and $V_{\text{out}}$ is a univalent vertex adjacent to a unique edge $E_{\text{out}}$. Then a {\em $d$-pointed marked tropical disk in $X_\Sigma$} is a map $h: \Gamma' \to M_\mathbb{R}$ satisfying the same conditions as above, except there is {\em no} balancing condition at $V_{\text{out}}$.

Now we fix $k$ points $P_1, P_2, \dots, P_k \in M_\mathbb{R} \cong \mathbb{R}^2$ in generic position. Then a {\em tropical disk in $(X_\Sigma; P_1, \dots, P_k)$ with boundary $Q$} is a $d$-pointed marked tropical disk in $X_\Sigma$ is a map $h: \Gamma' \to M_\mathbb{R}$ as above such that $h(V_{\text{out}}) = Q$ and $h(E_{p_j}) = P_{i_j}$ where $1 \leq i_1 < \dots < i_d \leq k$.
The {\em Maslox index} of such a disk $h$ is defined as $MI(h) := 2(N - d)$, where $N$ is the number of unbounded edges of $\Gamma'$ on which $h$ is non-constant (or number of unbounded edges in the image $h(\Gamma')$).

Taking the union of all Maslov index $0$ tropical disks, one obtains a scattering diagram $\mathscr{D} = \mathscr{D}(\Sigma; P_1, P_2, \dots, P_k)$, where $\Sigma$ denotes the fan of $\mathbb{P}^2$.
\begin{thm}[Proposition 4.7 in \cite{Gross10}]
	The scattering diagram $\mathscr{D}$ is consistent away from the points $P_1, P_2, \dots, P_k$, i.e. for any $P \in \text{Sing}(\mathscr{D})\setminus \{P_1, P_2, \dots, P_k\}$, we have $\Theta_\gamma = \text{Id}$ where $\Theta_\gamma$ is the path-ordered product along a loop $\gamma$ around $P$.
\end{thm}
Figure \ref{fig:scattering-MI0disks} below shows the cases $k = 1$ (left) and $k = 2$ (right).
\begin{figure}[ht]
	\includegraphics[scale=0.55]{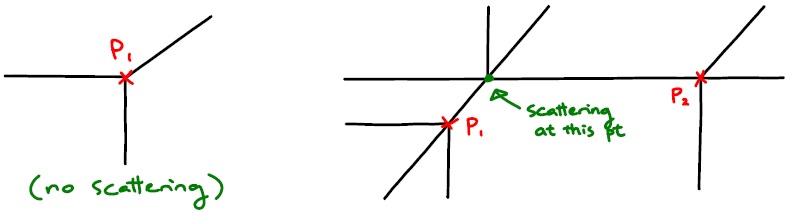}
	\caption{Scattering diagrams coming from loci of MI $0$ tropical disks.}\label{fig:scattering-MI0disks}
\end{figure}

Now the scattering diagram $\mathscr{D}$ divides $\mathbb{R}^2$ into different chambers (i.e. connected components of the complement $\mathbb{R}^2 \setminus \text{Supp}(\mathscr{D})$).
Then Gross' {\em $k$-pointed} perturbed superpotential is defined using counts of Maslov index $2$ tropical disks:
$$W_k(Q) := \sum_{\Gamma} Mono(\Gamma)$$
for $Q$ in a fixed chamber in $\mathbb{R}^2 \setminus \text{Supp}(\mathscr{D})$, where the sum is over all Maslov index 2 tropical disks passing through some of the points $P_1, P_2, \dots, P_k$ and with stop at $Q$, and $Mono(\Gamma)$ is a monomial in the variables $x, y$ associated to $\Gamma$. See \cite[Example 2.9]{Gross10} for some explicit examples.

When we move the stop $Q$ from one chamber to another, Gross proved that his perturbed superpotential displayed the following wall-crossing phenomenon:
\begin{thm}[Theorem 4.12 in \cite{Gross10}]
	For $Q_+, Q_- \in \mathbb{R}^2 \setminus \text{Supp}(\mathscr{D})$, we have
	$$W_k(Q_-) = \Theta_\gamma (W_k(Q_+)),$$
	where $\gamma$ is a path going from $Q_+$ to $Q_-$ and $\Theta_\gamma$ is the path-ordered product along $\gamma$, i.e. ordered product of wall-crossing factors associated to walls intersecting $\gamma$.
\end{thm}
\begin{figure}[ht]
	\includegraphics[scale=0.85]{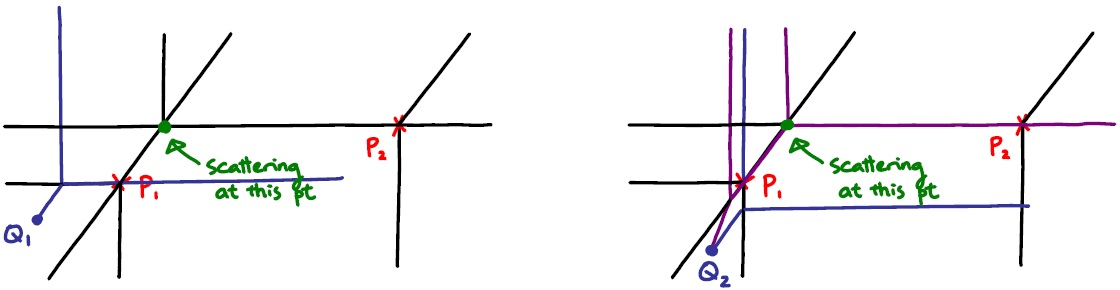}
	\caption{Wall-crossing for counts of MI $2$ tropical disks.}\label{fig:wall-crossing-disk-counts}
\end{figure}

To see how such tropical phenomena emerge from deformation theory, we consider the following DGBV algebra which controls (extended) deformations of the mirror Landau-Ginzburg model $(\check{X}, W)$:
$$(PV^{\bullet,\bullet}(\check{X}) := \Omega^{0,\bullet}(\check{X}, \wedge^\bullet T^{1,0}_{\check{X}}), \bar{\partial}_W := \bar{\partial} + [W, \cdot], \wedge, \Delta),$$
where the total degree on $PV^{\bullet,\bullet}(\check{X})$ is taken to be $\deg PV^{i,j} := j - i$, and $\Delta$ is the BV operator corresponding to $\partial$ on $\Omega^{\bullet,\bullet}(\check{X})$ under contraction by the standard holomorphic volume form $\check{\Omega} := d\log x \wedge d\log y$ on $\check{X} = (\mathbb{C}^*)^2$. Replacing $\bar{\partial}$ by the {\em twisted Dolbeault operator} $\bar{\partial}_W$, the Maurer-Cartan equation is written as
\begin{equation}\label{eqn:MC-eqn-LG}
	\bar{\partial}_W\Phi + \frac{1}{2}[\Phi, \Phi] = 0.
\end{equation}
In this case, Bogomolov-Tian-Todorov--type unobstructedness results have been proved by Katzarkov-Kontsevich-Pantev \cite{KKP17} (via degeneracy of a Hodge-to-de Rham spectral sequence as proved in Esnualt-Sabbah-Yu \cite{Esnault-Sabbah-Yu17}).

\begin{figure}[ht]
	\includegraphics[scale=0.6]{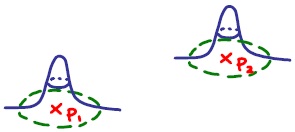}
	\caption{}\label{fig:input-LG}
\end{figure}

To solve the Maurer-Cartan equation \eqref{eqn:MC-eqn-LG}, we consider the input 
$$\Pi \in PV^{2,2}(\check{X}),$$
which is a sum of polyvector fields valued bump forms at each of the points $P_1, P_2, \dots, P_k$, as shown in Figure \ref{fig:input-LG}.

After choosing a suitable propagator $H$, we can use Kuranishi's method to solve \eqref{eqn:MC-eqn-LG} as in \S \ref{sec:scattering-from-MC}, namely, we can obtain a Maurer-Cartan solution by the formula:
\begin{equation*}\label{eqn:Kuranishi-formula-LG}
	\Phi := \Pi - H\left([W, \Phi] + \frac{1}{2} [\Phi, \Phi]\right)
\end{equation*}
(upon checking $p([W,\Phi] + \frac{1}{2}[\Phi, \Phi]) = 0$), which can in turn be expressed as a sum-over-trees formula $\Phi = \sum_{T:\text{trivalent trees}} W_T$, but now the decoration of each trivalent tree is modified as in Figure \ref{fig:trivalent-tree-LG}, namely, besides $\Pi$, we also align $W$ at some of the leaves of the tree $T$:
\begin{figure}[ht]
	\includegraphics[scale=0.55]{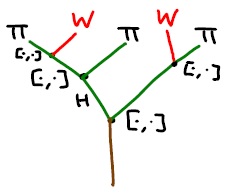}
	\caption{}\label{fig:trivalent-tree-LG}
\end{figure}

Now the interesting part is that such a Maurer-Cartan solution gives rise to all the tropical disks used by Gross:
\begin{thm}[\cite{Chan-Ma18}]\label{thm:Chan-Ma}
	Decomposing the Maurer-Cartan solution according to the double degrees: $\Phi = \Phi^{2,2} + \Phi^{1,1} + \Phi^{0,0}$, we have $\Phi^{2,2} = \Pi$, and
	\begin{equation*}
		\begin{split}
			\Phi^{1,1} & = \sum_{\text{$\Gamma$: MI $0$ tropical disks}} \alpha_\Gamma \log (\Theta_\Gamma),\\
			\Phi^{0,0} & = \sum_\text{$\Gamma$: MI $2$ tropical disks} \beta_\Gamma Mono(\Gamma),
		\end{split}
	\end{equation*}
	where $\alpha_\Gamma \in \Omega^{0,1}(\check{X})$ is a $(0,1)$-form supported near a wall in the scattering diagram $\mathscr{D} = \mathscr{D}(\Sigma; P_1, P_2, \dots, P_k)$, and $\beta_\Gamma \in \Omega^{0,0}(\check{X})$ is a function such that $\lim_{\hbar \to 0} \beta_\Gamma|_Q = 1$ for $Q$ in a fixed chamber in $\mathbb{R}^2 \setminus \mathscr{D}$.
\end{thm}
In particular, this gives the following one-to-one correspondences:
\begin{align*}
	\left\{\text{Maslov index $0$ tropical disks}\right\} & \longleftrightarrow \left\{\text{leading order terms in $\Phi^{1,1}$}\right\},\\
	\left\{\text{Maslov index $2$ tropical disks}\right\} & \longleftrightarrow \left\{\text{leading order terms in $\Phi^{0,0}$}\right\},
\end{align*}
where the latter is concerning the fixed chamber in $\mathbb{R}^2 \setminus \mathscr{D}$.

Away from the points $P_1, P_2, \dots, P_k$, the Maurer-Cartan equation \eqref{eqn:MC-eqn-LG} splits into the following two equations:
\begin{align}
	\bar{\partial} \Phi^{1,1} + \frac{1}{2}[\Phi^{1,1}, \Phi^{1,1}] & = 0, \label{eqn:MC-eqn-LG-I} \\
	\bar{\partial} \Phi^{0,0} + [\Phi^{1,1}, W + \Phi^{0,0}] & = 0. \label{eqn:MC-eqn-LG-II}
\end{align}
The results in \S \ref{sec:scattering-from-MC} say that the solution $\Phi^{1,1}$ of equation \eqref{eqn:MC-eqn-LG-I} gives rise to the consistent scattering diagram $\mathscr{D} = \mathscr{D}(\Sigma; P_1, P_2, \dots, P_k)$. Also, Theorem \ref{thm:Chan-Ma} says that, for $Q$ in a chamber of $\mathbb{R}^2 \setminus \mathscr{D}$, Gross' perturbed superpotential is given by $W_k(Q) = W(Q) + \Phi^{0,0}(Q)$. Now, the equation \eqref{eqn:MC-eqn-LG-II} can be rewritten as
\begin{equation}\label{eqn:MC-eqn-LG-III}
	(\bar{\partial} + \Phi^{1,1})(W + \Phi^{0,0}) = 0,
\end{equation}
and since $\Phi^{1,1}$ is gauge equivalent to $0$, i.e. there exists (uniquely) $\varphi \in \Omega^{0,0}(\check{X},T^{1,0}_{\check{X}})$ such that $\exp(\varphi) \ast 0 = \Phi^{1,1}$, the equation \eqref{eqn:MC-eqn-LG-III} is precisely telling us that
$$e^{\varphi}W_k(Q)$$
is a {\em global} holomorphic function (with respect to the original holomorphic structure on $\check{X}$) away from $P_1, P_2, \dots, P_k$. Furthermore, near a wall in $\mathscr{D}$, the gauge $\varphi$ is of the form (see Figure \ref{fig:wall-and-gauge} below):
$$\varphi = \left\{
\begin{array}{ll}
	\log(\Theta) & \text{on $H_{+}$,}\\
	0 & \text{on $H_{-}$.}
\end{array}\right.$$
Hence, this explains Gross' wall-crossing formula:
$$W_k(Q_-) = \Theta(W_k(Q_+))$$
across a single wall in $\mathscr{D}$.
\begin{figure}[ht]
	\includegraphics[scale=0.65]{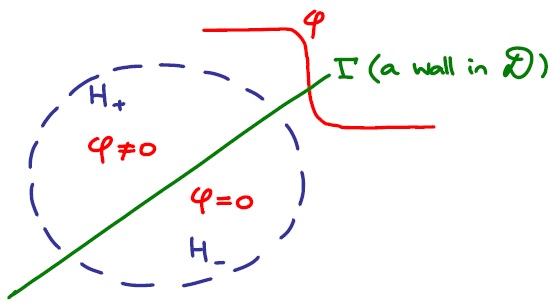}
	\caption{}\label{fig:wall-and-gauge}
\end{figure}

\section{Epilogue}

Without questions, the study of quantum corrections is of utmost importance in mirror symmetry.
Scattering diagrams and tropical counting are descriptions of quantum corrections using the language in the tropical world. So no doubt they play important roles in mirror symmetry, and in particular, in the reconstruction problem. Recently, such tropical objects have even found applications beyond mirror symmetry -- for example, in the study of cluster algebras and cluster varieties; see e.g. \cite{GHK13, GHKK18}.

On the other hand, the lesson we learn from the results described in this note is that consistent scattering diagrams and tropical counts are actually encoded in Maurer-Cartan solutions of the Kodaira-Spencer DGLA. Thus, for the purpose of proving mirror statements or even understanding mirror symmetry, one may bypass the complicated combinatorics of these objects. This was the motivation behind the recent works \cite{Chan-Leung-Ma19, Chan-Ma19, Chan-Ma-Suen20}, where we directly constructed the DGLA (or better, DGBV algebra) that governs smoothing of degenerate Calabi-Yau varieties and prove existence of smoothing without using scattering diagrams.

Besides giving a unified proof and substantial extension of previous smoothing results \cite{Felten-Filip-Ruddat19}, the algebraic framework developed in \cite{Chan-Leung-Ma19} can actually be combined with the techniques of asymptotic analysis in \cite{Chan-Leung-Ma17, Chan-Ma18} to globalize the results reviewed in this note, which are local in nature. In the forthcoming work \cite{Chan-Leung-Ma21}, we will prove the original, global version of Fukaya's conjecture \cite[Conjecture 5.3]{Fukaya05} (reformulated as discussed in this note) using this strategy. To speculate further, these methods will yield a local-to-global approach for proving genus $0$ mirror symmetry. We hope to report on this in the near future.

\section*{Acknowledgment}
The author would like to thank Naichung Conan Leung and Ziming Nikolas Ma for very fruitful collaborations which lead to much of the work described here, and the referee for valuable comments and pointing out a number of errors.

Research described in this paper were substantially supported by grants from the Research Grants Council of the Hong Kong Special Administrative Region, China (Project No. CUHK14303019 \& CUHK14301420).

\bibliographystyle{amsplain}
\bibliography{geometry}

\providecommand{\bysame}{\leavevmode\hbox to3em{\hrulefill}\thinspace}
\providecommand{\MR}{\relax\ifhmode\unskip\space\fi MR }
\providecommand{\MRhref}[2]{%
  \href{http://www.ams.org/mathscinet-getitem?mr=#1}{#2}
}
\providecommand{\href}[2]{#2}
\begin{thebibliography}{10}

\bibitem{Batyrev93}
V.~Batyrev, \emph{Quantum cohomology rings of toric manifolds}, Ast\'erisque
  (1993), no.~218, 9--34, Journ{\'e}es de G{\'e}om{\'e}trie Alg{\'e}brique
  d'Orsay (Orsay, 1992).

\bibitem{Bogomolov78}
F.~A. Bogomolov, \emph{Hamiltonian {K}\"{a}hlerian manifolds}, Dokl. Akad. Nauk
  SSSR \textbf{243} (1978), no.~5, 1101--1104.

\bibitem{Chan-Leung10a}
K.~Chan and N.~C. Leung, \emph{Mirror symmetry for toric {F}ano manifolds via
  {SYZ} transformations}, Adv. Math. \textbf{223} (2010), no.~3, 797--839.

\bibitem{Chan-Leung-Ma19}
K.~Chan, N.~C. Leung, and Z.~N. Ma, \emph{Geometry of the {M}aurer-{C}artan
  equation near degenerate {C}alabi-{Y}au varieties}, J. Differential Geom., to
  appear, \href{http://arxiv.org/abs/1902.11174}{arXiv:1902.11174}.

\bibitem{Chan-Leung-Ma21}
\bysame, \emph{Smoothing, scattering, and a conjecture of {F}ukaya}, preprint
  (2022).

\bibitem{Chan-Leung-Ma17}
\bysame, \emph{Scattering diagrams from asymptotic analysis on
  {M}aurer-{C}artan equations}, J. Eur. Math. Soc. (JEMS) \textbf{24} (2022),
  no.~3, 773--849.

\bibitem{Chan-Ma18}
K.~Chan and Z.~N. Ma, \emph{Tropical counting from asymptotic analysis on
  {M}aurer-{C}artan equations}, Trans. Amer. Math. Soc. \textbf{373} (2020),
  no.~9, 6411--6450.

\bibitem{Chan-Ma19}
\bysame, \emph{Smoothing pairs over degenerate {C}alabi-{Y}au varieties}, Int.
  Math. Res. Not. IMRN (2022), no.~4, 2582--2614.

\bibitem{Chan-Ma-Suen20}
K.~Chan, Z.~N. Ma, and Y.-H. Suen, \emph{Tropical {L}agrangian multi-sections
  and smoothing of locally free sheaves over degenerate {C}alabi-{Y}au
  surfaces}, Adv. Math. \textbf{401} (2022), Paper No. 108280, 37pp.

\bibitem{Cho-Oh06}
C.-H. Cho and Y.-G. Oh, \emph{Floer cohomology and disc instantons of
  {L}agrangian torus fibers in {F}ano toric manifolds}, Asian J. Math.
  \textbf{10} (2006), no.~4, 773--814.

\bibitem{Duistermaat80}
J.~J. Duistermaat, \emph{On global action-angle coordinates}, Comm. Pure Appl.
  Math. \textbf{33} (1980), no.~6, 687--706.

\bibitem{Esnault-Sabbah-Yu17}
H.~Esnault, C.~Sabbah, and J.-D. Yu, \emph{{$E_1$}-degeneration of the
  irregular {H}odge filtration}, J. Reine Angew. Math. \textbf{729} (2017),
  171--227, With an appendix by Morihiko Saito.

\bibitem{Felten-Filip-Ruddat19}
S.~Felten, M.~Filip, and H.~Ruddat, \emph{Smoothing toroidal crossing spaces},
  Forum Math. Pi \textbf{9} (2021), Paper No. e7, 36.

\bibitem{Fukaya05}
K.~Fukaya, \emph{Multivalued {M}orse theory, asymptotic analysis and mirror
  symmetry}, Graphs and patterns in mathematics and theoretical physics, Proc.
  Sympos. Pure Math., vol.~73, Amer. Math. Soc., Providence, RI, 2005,
  pp.~205--278.

\bibitem{FOOO-book}
K.~Fukaya, Y.-G. Oh, H.~Ohta, and K.~Ono, \emph{Lagrangian intersection {F}loer
  theory: anomaly and obstruction.}, AMS/IP Studies in Advanced Mathematics,
  vol.~46, American Mathematical Society, Providence, RI, 2009.

\bibitem{FOOO-toricI}
\bysame, \emph{Lagrangian {F}loer theory on compact toric manifolds. {I}}, Duke
  Math. J. \textbf{151} (2010), no.~1, 23--174.

\bibitem{FOOO-toricII}
\bysame, \emph{Lagrangian {F}loer theory on compact toric manifolds {II}: bulk
  deformations}, Selecta Math. (N.S.) \textbf{17} (2011), no.~3, 609--711.

\bibitem{FOOO-toricIII}
\bysame, \emph{Lagrangian {F}loer theory and mirror symmetry on compact toric
  manifolds}, Ast\'erisque (2016), no.~376, vi+340.

\bibitem{Givental98}
A.~Givental, \emph{A mirror theorem for toric complete intersections},
  Topological field theory, primitive forms and related topics ({K}yoto, 1996),
  Progr. Math., vol. 160, Birkh\"auser Boston, Boston, MA, 1998, pp.~141--175.

\bibitem{Gross10}
M.~Gross, \emph{Mirror symmetry for {$\Bbb P^2$} and tropical geometry}, Adv.
  Math. \textbf{224} (2010), no.~1, 169--245.

\bibitem{GHK13}
M.~Gross, P.~Hacking, and S.~Keel, \emph{Birational geometry of cluster
  algebras}, Algebr. Geom. \textbf{2} (2015), no.~2, 137--175.

\bibitem{GHKK18}
M.~Gross, P.~Hacking, S.~Keel, and M.~Kontsevich, \emph{Canonical bases for
  cluster algebras}, J. Amer. Math. Soc. \textbf{31} (2018), no.~2, 497--608.

\bibitem{GPS10}
M.~Gross, R.~Pandharipande, and B.~Siebert, \emph{The tropical vertex}, Duke
  Math. J. \textbf{153} (2010), no.~2, 297--362.

\bibitem{Gross-Siebert03}
M.~Gross and B.~Siebert, \emph{Affine manifolds, log structures, and mirror
  symmetry}, Turkish J. Math. \textbf{27} (2003), no.~1, 33--60.

\bibitem{Gross-Siebert-logI}
\bysame, \emph{Mirror symmetry via logarithmic degeneration data. {I}}, J.
  Differential Geom. \textbf{72} (2006), no.~2, 169--338.

\bibitem{Gross-Siebert-logII}
\bysame, \emph{Mirror symmetry via logarithmic degeneration data, {II}}, J.
  Algebraic Geom. \textbf{19} (2010), no.~4, 679--780.

\bibitem{Gross-Siebert-reconstruction}
\bysame, \emph{From real affine geometry to complex geometry}, Ann. of Math.
  (2) \textbf{174} (2011), no.~3, 1301--1428.

\bibitem{Hori-Vafa00}
K.~Hori and C.~Vafa, \emph{Mirror symmetry}, preprint (2000),
  \href{http://arxiv.org/abs/hep-th/0002222}{arXiv:hep-th/0002222}.

\bibitem{KKP17}
L.~Katzarkov, M.~Kontsevich, and T.~Pantev, \emph{Bogomolov-{T}ian-{T}odorov
  theorems for {L}andau-{G}inzburg models}, J. Differential Geom. \textbf{105}
  (2017), no.~1, 55--117.

\bibitem{Kodaira-Spencer-I-II}
K.~Kodaira and D.~C. Spencer, \emph{On deformations of complex analytic
  structures. {I}, {II}}, Ann. of Math. (2) \textbf{67} (1958), 328--466.

\bibitem{Kodaira-Spencer-III}
\bysame, \emph{On deformations of complex analytic structures. {III}.
  {S}tability theorems for complex structures}, Ann. of Math. (2) \textbf{71}
  (1960), 43--76.

\bibitem{Kontsevich-Soibelman01}
M.~Kontsevich and Y.~Soibelman, \emph{Homological mirror symmetry and torus
  fibrations}, Symplectic geometry and mirror symmetry ({S}eoul, 2000), World
  Sci. Publ., River Edge, NJ, 2001, pp.~203--263.

\bibitem{Kontsevich-Soibelman06}
\bysame, \emph{Affine structures and non-{A}rchimedean analytic spaces}, The
  unity of mathematics, Progr. Math., vol. 244, Birkh\"auser Boston, Boston,
  MA, 2006, pp.~321--385.

\bibitem{Kuranishi65}
M.~Kuranishi, \emph{New proof for the existence of locally complete families of
  complex structures}, Proc. {C}onf. {C}omplex {A}nalysis ({M}inneapolis,
  1964), Springer, Berlin, 1965, pp.~142--154.

\bibitem{Mikhalkin05}
G.~Mikhalkin, \emph{Enumerative tropical algebraic geometry in {$\Bbb R^2$}},
  J. Amer. Math. Soc. \textbf{18} (2005), no.~2, 313--377.

\bibitem{Nishinou12}
T.~Nishinou, \emph{Disk counting on toric varieties via tropical curves}, Amer.
  J. Math. \textbf{134} (2012), no.~6, 1423--1472.

\bibitem{Nishinou-Siebert06}
T.~Nishinou and B.~Siebert, \emph{Toric degenerations of toric varieties and
  tropical curves}, Duke Math. J. \textbf{135} (2006), no.~1, 1--51.

\bibitem{SYZ96}
A.~Strominger, S.-T. Yau, and E.~Zaslow, \emph{Mirror symmetry is
  {$T$}-duality}, Nuclear Phys. B \textbf{479} (1996), no.~1-2, 243--259.

\bibitem{Tian86}
G.~Tian, \emph{Smoothness of the universal deformation space of compact
  {C}alabi-{Y}au manifolds and its {P}etersson-{W}eil metric}, Mathematical
  aspects of string theory ({S}an {D}iego, {C}alif., 1986), Adv. Ser. Math.
  Phys., vol.~1, World Sci. Publishing, Singapore, 1987, pp.~629--646.

\bibitem{Todorov89}
A.~N. Todorov, \emph{The {W}eil-{P}etersson geometry of the moduli space of
  {${\rm SU}(n\geq 3)$} ({C}alabi-{Y}au) manifolds. {I}}, Comm. Math. Phys.
  \textbf{126} (1989), no.~2, 325--346.

\end{thebibliography}

\end{document}